\input amstex
\input amsppt.sty

\define\x{\frak X}
\define\pr{\operatorname{Prob}}
\define\fr{\operatorname{Fr}}
\define\C{\Bbb C}
\define\R{\Bbb R}
\define\Z{\Bbb Z}
\define\N{\Bbb N}
\define\tht{\thetag}
\define\la{\lambda}
\define\ze{\zeta}
\define\r{\operatorname{Res}\limits_{\zeta=x}}
\define\wt{\widetilde}

\NoRunningHeads
\TagsOnRight
\magnification 1200
\hoffset=0.5cm
\voffset=-1.5cm

\topmatter
\title 
Asymptotic representation theory and Riemann--Hilbert problem
\endtitle
\author Alexei Borodin
\endauthor

\abstract We show how the Riemann--Hilbert problem can be used to
compute correlation kernels for determinantal point processes arising in different models of asymptotic combinatorics and representation theory.
The Whittaker kernel and the discrete Bessel kernel are computed as examples. \endabstract

\endtopmatter

\document

\head Introduction
\endhead

A (discrete or continuous) random point process is called {\it determinantal}
if its correlation functions have the form
$$
\rho_n(x_1,\dots,x_n)=\det[K(x_i,x_j)]_{i,j=1}^n,
$$
where $K(x,y)$ is a function in two variables called the {\it correlation
kernel}. A major source of such point processes is Random Matrix Theory.
All the ``unitary'' or ``$\beta=2$'' ensembles of random matrices lead to
determinantal point processes which describe the eigenvalues of these
matrices.

Determinantal point processes
also arise naturally in problems of asymptotic combinatorics and asymptotic representation theory, see \cite{BO1}--\cite{BO4}, \cite{BOO}, \cite{J}, \cite{Ol2}. Usually, it is not very hard to see that the process that we are
interested in is determinantal. A harder problem is to compute the correlation
kernel of this process explicitly. The goal of this paper is to give an
informal introduction to a new method of obtaining explicit formulas for correlation kernels.
It should be emphasized that in representation theoretic models which we
consider the kernels cannot be expressed through orthogonal polynomials, as it
often happens in random matrix models. That is why we had to invent something different.

The heart of the method is the {\it Riemann--Hilbert problem} (RHP, for
short). This is a classical problem which consists of factorizing a
matrix--valued function on a contour in the complex plane into a product of a
function which is holomorphic inside the contour and a function which is
holomorphic outside the contour. It turns out that the problem of computing the
correlation kernels can be reduced to solving a RHP of a rather special form.
The input of the RHP (the function to be factorized) is always rather simple
and can be read off the representation theoretic quantities such as dimensions
of irreducible representations of the corresponding groups.
We also employ a discrete analog of RHP described in \cite{B2}. 

The special form of our concrete RHPs allows us to reduce them to certain linear
ordinary differential equations (this is the key step), which have classical
special functions as their solutions. This immediately leads to explicit
formulas for the needed correlation kernels.

The approach also happens to be very effective for the derivation of
(nonlinear ordinary differential) Painlev\'e equations describing the ``gap
probabilities'' in both random matrix and representation theoretic models, see \cite{BD}, \cite{B3}.
However, this issue will not be addressed in this paper.

The paper is organized as follows. In \S1 we explain what a determinantal point process is and give a couple of examples. In \S2 we argue that in many models correlation kernels give rise to what is called ``integrable integral operators''. In \S3 we relate integrable operators to RHP. In \S4 we derive the Whittaker kernel arising in a problem of harmonic analysis on the infinite symmetric group. In \S5 we derive the discrete Bessel kernel associated with the poissonized Plancherel measures on symmetric groups.

This paper is an expanded text of lectures the author gave at the NATO Advanced Study Institute ``Asymptotic combinatorics with applications to mathematical physics'' in July 2001 in St. Petersburg.
It is a great pleasure to thank the organizers for the invitation and
for the warm hospitality. The author would also like to thank Grigori Olshanski and Percy Deift for helpful discussions.

This research was partially conducted during the period the author served as a Clay Mathematics Institute Long-Term Prize Fellow. This work was also partially supported by the NSF grant DMS-9729992. 
 
\head 1. Determinantal point processes
\endhead

\example{Definition 1.1} Let $\x$ be a discrete space. A probability measure on $2^\x$ is called a {\it determinantal point process} if there exists a function $K:\x\times \x\to\C$ such that
$$
\pr\left\{A\in 2^\x\,|\,A\supset \{x_1,\dots,x_n\}\right\}
=\det[K(x_i,x_j)]_{i,j=1}^n
$$
for any finite subset $\{x_1,\dots.x_n\}$ of $\x$. The function $K$ is called the {\it correlation kernel}. The functions
$$
\gathered
\rho_n: \{\text{$n$-point subsets of $\x$}\}\to [0,1]\\
\rho_n:\{x_1,\dots,x_n\}\mapsto \pr\{A\,|\,A\supset \{x_1,\dots,x_n\}\}
\endgathered
$$
are called the {\it correlation functions}.
\endexample

\example{Example 1.2} Consider a kernel $L:\x \times \x\to \C$
such that
\smallskip
$\bullet$ $\det[L(x_i,x_j)]_{i,j=1}^k\ge 0$ for all $k$-point subsets $\{y_1,\dots,y_k\}$ of $\x$. 
\smallskip
$\bullet$ $L$ defines a trace class operator in $\ell^2(\x)$, for example, $\sum_{x,y\in\x}|L(x,y)|<\infty$ or $L$ is finite rank. In particular, this condition is empty if $|\x|<\infty$.
\medskip
Set
$$
\pr\left\{\{y_1,\dots,y_k\}\right\}=\frac{1}{\det(1+L)}\cdot\det[L(y_i,y_j)]_{i,j=1}^k\,.
$$ 
This defines a probability measure on $2^\x$ concentrated on finite subsets. Moreover, this defines a determinantal point process. The correlation kernel $K(x,y)$ is equal to the matrix of the operator $K=L(1+L)^{-1}$ acting on $\ell^2(\x)$. See \cite{DVJ}, \cite{BOO, Appendix} for details.  
\endexample

\example{Definition 1.3} Let $\x$ be a finite or infinite interval inside
$\R$ (e.g., $\R$ itself). A probability measure on locally finite subsets of $\x$ is called a {\it determinantal point process} if there exists a function $K:\x\times\x\to \C$ such that
$$
\multline
\lim_{\Delta x_1,\dots,\Delta x_n\to 0}
\frac{\pr\{A\in 2^\x_{\text{loc.fin.}}\,|\,
\text{$A$ intersects $[x_i,x_i+\Delta x_i]$ for all $i=1,\dots,n$}\}}{\Delta x_1\cdots\Delta x_n}\\
=\det[K(x_i,x_j)]_{i,j=1}^n
\endmultline
$$
for any finite subset $\{x_1,\dots,x_n\}$ of $\x$. The function $K$ is called the {\it correlation kernel} and the left--hand side of the equality above is called the {\it $n$th correlation function}.
\endexample
\example{Example 1.4} Let $w(x)$ be a positive function on $\x$ such that
all the moments $\int_\x x^nw(x)dx$ are finite. Pick a number $N\in \N$
and define a probability measure on $N$-point subsets of $\x$ by the formula 
$$
P_N(dx_1,\dots,dx_N)=c_N \prod_{1\le i<j\le N}(x_i-x_j)^2\prod_{1\le k\le N} w(x_k)dx_k.
$$
Here $c_N>0$ is a normalizing constant. This is a determinantal point process. The correlation kernel is equal to the $N$th Christoffel--Darboux kernel $K_N(x,y)$ associated with $w(x)$, multiplied by $\sqrt{w(x)w(y)}$. That is, let 
$$
p_0=1,\ p_1(x),\ p_2(x),\dots
$$
be monic (= leading coefficient 1) orthogonal polynomials on $\x$ with the weight function $w(x)$:
$$
\gathered
p_m(x)=x^m+\text{ lower order terms },\\
\int_\x p_m(x) p_n(x)w(x)dx=h_m\delta_{mn}, \quad m,n=0,1,2,\dots\,.
\endgathered
$$
Then the correlation kernel is equal to 
$$
\multline
K_N(x,y)=\sum_{k=0}^N\frac{p_k(x)p_k(y)}{h_k}\,\sqrt{w(x)w(y)}\\
=\frac 1{h_{N-1}}\,\frac{p_N(x)p_{N-1}(y)-p_{N-1}(x)p_N(y)}{x-y}\,\sqrt{w(x)w(y)}.
\endmultline
$$
The construction of this example also makes sense in the discrete setting. See \cite{Dy}, 
\cite{Me}, \cite{NW}, \cite{J} for details. 
\endexample

\example{Remark 1.5} The correlation kernel of a determinantal point process is not defined uniquely! In particular, transformations of the form $K(x,y)\to \frac{f(x)}{f(y)}K(x,y)$  do not change the correlation functions.
\endexample

\head 2. Correlation kernels as integrable operators
\endhead
Observe that the kernel $K_N(x,y)$ of Example 1.4 has the form 
$$
K_N(x,y)=\frac{\phi(x)\psi(y)-\psi(x)\phi(y)}{x-y}
$$
for appropriate $\phi$ and $\psi$. Most kernels appearing in  ``$\beta=2$ ensembles'' of Random Matrix Theory have this form, because they are either kernels of Christoffel--Darboux type as in Example 1.4 above, or scaling limits of such kernels. However, it is an experimental fact that integral operators with such kernels appear in many different areas of mathematics, see \cite{De}. 

\example{Definition 2.1}
An integral operator with kernel of the form
$$
\frac{f_1(x)g_1(y)+f_2(x)g_2(y)}{x-y}
\tag 2.1
$$
is called {\it integrable}. Here we assume that $f_1(x)g_1(x)+f_2(x)g_2(x)=0$ so that there is no singularity on the diagonal. Diagonal values of the kernel are then defined by continuity.

The class of integrable operators was singled out in the work of Its, Izergin, Korepin, and Slavnov on quantum inverse scattering method in 1990 \cite{IIKS}.

We will also call an operator acting in the $\ell^2$-space on a discrete space {\it integrable} if its matrix has the form \tht{2.1}. It is not obvious how to define the diagonal entries of a discrete integrable operator in general. However, in all concrete situations we are aware of, this question has a natural answer.
\endexample
\example{Example 2.2 (poissonized Plancherel measure, cf. \cite{BOO})} Consider the probability measure on the set of all Young diagrams given by the formula
$$
\pr\{\la\}=e^{-\theta}\theta^{|\la|}{\left(\frac{\dim\la}{|\la|!}\right)}^2.
\tag 2.2
$$
Here $\theta>0$ is a parameter, $\dim\la$ is the number of standard Young tableaux of shape $\la$ or the dimension of the irreducible representation of the symmetric group $S_{|\la|}$ corresponding to $\la$.
Denote by $(p_1,\dots,p_d|q_1,\dots,q_d)$ the Frobenius coordinates of $\la$ (see \cite{Ma, \S1} for the definition of Frobenius coordinates).
Here $d$ is the number of diagonal boxes in $\la$.
Set $\Z'=\Z+\frac 12=\{\pm\frac 12,\pm\frac 32,\dots\}$. 

Let us associate to any Young diagram $\la=(p\,|\,q)$ a point configuration $\fr(\la)\subset \Z'$ as follows:
$$
\fr(\la)=\left\{p_1+\tfrac 12,\dots,p_d+\tfrac 12,-q_1-\tfrac 12,\dots,-q_d-\tfrac 12\right\}.
$$

It turns out that together with \tht{2.2} this defines a determinantal point process on $\Z'$. Indeed, the well-known hook formula for $\dim\la$  easily implies
$$
\pr\{\la\}=e^{-\theta}\left(\det\left[\frac{\theta^{\frac{p_i+q_j}2}}{(p_i-\frac 12)!(q_j-\frac 12)!(p_i+q_j)}\right]_{i,j=1}^d\right)^2=e^{-\theta}
\det[L(y_i,y_j)]_{i,j=1}^{2d}
$$
where $\{y_1,\dots,y_{2d}\}=\fr(\la)$, and $L(x,y)$ is a $\Z'\times\Z'$ matrix defined by
$$
L(x,y)=\cases 0,&\text{if  } xy>0,\\
\dfrac{\theta^{\frac{|x|+|y|}2}}{(|x|-\frac 12)!(|y|-\frac 12)!}\,\dfrac 1{x-y}\,,&
\text{if  }xy<0.
\endcases
$$
In the block form corresponding to the splitting $\Z'=\Z'_+\sqcup \Z'_-$ it looks as follows
$$
L(x,y)=\bmatrix 0&\dfrac{\theta^{\frac{x-y}2}}{(x-\frac 12)!(-y-\frac 12)!}\,\dfrac 1{x-y}\\
\dfrac{\theta^{\frac{-x+y}2}}{(-x-\frac 12)!(y-\frac 12)!}\,\dfrac 1{x-y}&0
\endbmatrix.
$$
The kernel $L(x,y)$ belongs to the class of integrable kernels. Indeed, if we set
$$
f_1(x)=g_2(y)=\cases \dfrac{\theta^{\frac x2}}{(x-\frac 12)!},&x>0,\\0,&x<0,
\endcases
\qquad 
f_2(x)=g_1(y)=\cases 0,&x>0,\\
\dfrac{\theta^{-\frac x2}}{(-x-\frac 12)!},&x<0,
\endcases
$$
then it is immediately verified that $L(x,y)=(f_1(x)g_1(y)+f_2(x)g_2(y))/(x-y)$.
Comparing the formulas with Example 1.2, we also conclude that $e^\theta=\det(1+L)$.\footnote{Since $\sum_{x,y\in\Z'}|L(x,y)|<\infty$, the operator $L$ is trace class, and $\det(1+L)$ is well--defined.}
\endexample

What we see in this example is that $L$ is an integrable kernel. We also 
know, see Example 1.2, that the correlation kernel $K$ is given by $K=L(1+L)^{-1}$. Is this kernel also integrable? The answer is positive; the general claim in the continuous case was proved in \cite{IIKS}, the discrete case was worked out in \cite{B2}. 

Furthermore, it turns out that in many situations there is an algorithm of computing the correlation kernel $K$ if $L$ is an integrable kernel which is ``simple enough''. The algorithm is based on a classical problem of complex analysis called the {\it Riemann--Hilbert problem} (RHP, for short).
 
Let us point out that our algorithm is not applicable to deriving correlation kernels in the ``$\beta=2$'' model of Random Matrix Theory. Indeed, the Christoffel--Darboux kernels have norm 1, since they are just projection operators. Thus, it is impossible to define the kernel $L=K(1-K)^{-1}$, because $(1-K)$ is not invertible. In this sense, RMT deals with ``degenerate'' determinantal point processes.

On the other hand, the orthogonal polynomial method of computing the correlation kernels, which has been so successful in RMT, cannot be applied directly to the representation theoretic models like Example 2.2 above (see, however, \cite{J}). The algorithm explained below may be viewed as a substitute for this method.

\head 3. Riemann--Hilbert problem
\endhead
Let $\Sigma$ be an oriented contour in $\C$. We agree that $(+)$-side is on the left of the contour, and $(-)$-side is on the right of the contour. Let $v$ be a $2\times 2$-matrix valued function on $\Sigma$. 

\example{Definition 3.1} We say that a matrix function $m:\C\setminus \Sigma\to\operatorname{Mat}(2,\C)$ solves the RHP $(\Sigma,v)$ if 

(1) $m$ is analytic in $\C\setminus \Sigma$;

(2) $m_+=m_-v$ on $\Sigma$, where $m_{\pm}(x)=\lim\limits_{\zeta\to x
\text{  from $(\pm)$-side}} m(\zeta)$.

We say that $m$ solves the {\it normalized} RHP $(\Sigma,v)$ if, in addition, we have

(3) $m(\zeta)\to I=\bmatrix 1&0\\0& 1\endbmatrix $ as $\ze\to\infty$.
\endexample

Next we explain what is a {\it discrete} Riemann--Hilbert problem (DRHP, for short).

Let $X$ be a locally finite subset of $\C$, and let $w$ be a $2\times 2$-matrix valued function on $X$. 

\example{Definition 3.2} 
We say that a matrix function $m:\C\setminus X\to\operatorname{Mat}(2,C)$ solves the DRHP $(X,w)$ if

(1) $m$ is analytic in $\C\setminus X$;

(2) $m$ has simple poles at the points of $X$, and
$$
\r m(\zeta)=\lim_{\ze\to x}(m(\zeta)w(x))\quad\text{  for any  } x\in X.
$$

We say that $m$ solves the {\it normalized} DRHP $(X,w)$ if

(3) $m(\zeta)\to I=\bmatrix 1&0\\0& 1\endbmatrix $ as $\ze\to\infty$.

If the set $X$ is infinite, the last relation should hold when the distance from $\zeta$ to $X$ is bounded away from zero. 
\endexample

Our next step is to explain how to reduce, for an integrable operator $L$, the computation of the operator $K=L(1+L)^{-1}$ to a (discrete or continuous) RHP.

\subhead Continuous picture \cite{IIKS}\endsubhead
Let $L$ be an integrable operator on $L^2(\Sigma, |d\ze|)$, $\Sigma\subset \C$, with the kernel ($x,y\in\Sigma$)
$$
L(x,y)=\frac{f_1(x)g_1(y)+f_2(x)g_2(y)}{x-y}\,,\, \qquad 
f_1(x)g_1(x)+f_2(x)g_2(x)\equiv 0.
$$
Assume that $(1+L)$ is invertible. 
\proclaim{Theorem 3.3} There exists a unique solution of the normalized RHP $(\Sigma,v)$ with 
$$
v=I+2\pi i \bmatrix f_1\\f_2\endbmatrix \bmatrix g_1 \ g_2\endbmatrix=\bmatrix 1+2\pi i f_1g_1& 2\pi i f_1g_2\\
2\pi i f_2g_1& 1+2\pi i f_2g_2\endbmatrix.
$$
For $x\in \Sigma$ set 
$$
\gathered
\bmatrix F_1(x)\\F_2(x)\endbmatrix=\lim_{\ze\to x}m(\ze)
\bmatrix f_1(x)\\f_2(x)\endbmatrix,\\
\bmatrix G_1(x)\\G_2(x)\endbmatrix=\lim_{\ze\to x}m^{-t}(\ze)
\bmatrix g_1(x)\\g_2(x)\endbmatrix.
\endgathered
$$
Then the kernel of the operator $K=L(1+L)^{-1}$ has the form $(x,y\in\Sigma)$
$$
K(x,y)=\frac{F_1(x)G_1(y)+F_2(x)G_2(y)}{x-y}\, \quad \text{  and  }\quad
F_1(x)G_1(x)+F_2(x)G_2(x)\equiv 0.
$$
\endproclaim
\example{Example 3.4} Let $\Sigma$ be a simple closed curve in $\C$ oriented clockwise (so that the $(+)$-side is outside $\Sigma$), and let $L$ be an integrable operator such that the functions $f_1,f_2,g_1,g_2$ can be extended to analytic functions inside $\Sigma$.
Then the solution of the normalized RHP $(\Sigma,v)$ has the form
$$
m=\cases \bmatrix 1&0\\0&1\endbmatrix &\text{outside $\Sigma$},\\
I-2\pi i \bmatrix f_1\\f_2\endbmatrix \bmatrix g_1 \ g_2\endbmatrix &\text{inside $\Sigma$}.
\endcases
$$
Then we immediately obtain $F_i=f_i$, $G_i=g_i$, $i=1,2$; and $K=L(1+L)^{-1}=L$. On the other hand, this is obvious because
$\int_\Sigma L(x,y)L(y,z)dy=0$ by Cauchy's theorem which means that $L^2=0$.
\endexample

\subhead Discrete picture \cite{B2} \endsubhead
Let $L$ be an integrable operator on $\ell^2(X)$, $X\subset \C$, with the kernel
$$
L(x,y)=\cases\dfrac{f_1(x)g_1(y)+f_2(x)g_2(y)}{x-y},&x\ne y,\\
0&x=y,\endcases
$$
with $f_1(x)g_1(x)+f_2(x)g_2(x)\equiv 0$. Assume that $(1+L)$ is invertible.
\proclaim{Theorem 3.5} There exists a unique solution of the normalized DRHP $(X,w)$ with 
$$
w=-\bmatrix f_1\\f_2\endbmatrix \bmatrix g_1 \ g_2\endbmatrix=\bmatrix - f_1g_1& - f_1g_2\\
- f_2g_1& -f_2g_2\endbmatrix.
$$
For $x\in X$ set 
$$
\gathered
\bmatrix F_1(x)\\F_2(x)\endbmatrix=\lim_{\ze\to x}m(\ze)
\bmatrix f_1(x)\\f_2(x)\endbmatrix,\\
\bmatrix G_1(x)\\G_2(x)\endbmatrix=\lim_{\ze\to x}m^{-t}(\ze)
\bmatrix g_1(x)\\g_2(x)\endbmatrix.
\endgathered
$$
Then the kernel of the operator $K=L(1+L)^{-1}$ has the form $(x,y\in\Sigma)$
$$
K(x,y)=\cases \dfrac{F_1(x)G_1(y)+F_2(x)G_2(y)}{x-y},&x\ne y,\\
\bmatrix G_1(x) \ G_2(x)\endbmatrix \lim\limits_{\ze\to x}\left(m'(\zeta)\bmatrix f_1(x)\\f_2(x)\endbmatrix\right),&x=y.
\endcases
$$
We also have $F_1(x)G_1(x)+F_2(x)G_2(x)\equiv 0$ on $X$.
\endproclaim
Theorem 3.5 can be extended to the case when $L(x,x)\ne 0$, see \cite{B2, Remark 4.2}.

\example{Example 3.6} Let $X=\{a,b\}$ be a two-point subset of $\C$, and 
$$
L=\bmatrix 0 & \mu\\ \nu&0\endbmatrix.
$$
Then $L$ is integrable with 
$$
f_1=\cases 0\\ \nu(b-a)\endcases,\quad f_2=\cases \mu(a-b)\\0\endcases,\quad g_1=\cases 1\\ 0\endcases,\quad g_2=\cases 0\\1\endcases.
$$
The notation means that, say, $f_1(a)=0,\ f_1(b)=\nu(b-a)$. Then 
$$
w(a)=\bmatrix 0&0\\ \mu(b-a)&0\endbmatrix,\quad
w(b)=\bmatrix 0&\nu (a-b)\\ 0&0\endbmatrix.
$$
Then the matrix $m(\ze)$ has the form
$$
m(\ze)=I+\frac 1{\ze-a}\,\frac{\mu(a-b)}{1-\mu\nu}\bmatrix \nu&0\\-1&0\endbmatrix+\frac 1{\ze-b}\,\frac{\nu(b-a)}{1-\mu\nu}
\bmatrix 0&-1\\0&\mu\endbmatrix.
$$
One can check that $\det m\equiv 1$, and 
$$
m^{-t}(\ze)=I+\frac 1{\ze-a}\,\frac{\mu(a-b)}{1-\mu\nu}\bmatrix 0&1\\0&\nu\endbmatrix+\frac 1{\ze-b}\,\frac{\nu(b-a)}{1-\mu\nu}
\bmatrix \mu&0\\1&0\endbmatrix.
$$
Further, 
$$
\gathered
F_1=\cases \frac{\mu\nu(b-a)}{1-\mu\nu}\\ \frac{\nu(b-a)}{1-\mu\nu}\endcases,\quad F_2=\cases\frac{ \mu(a-b)}{1-\mu\nu} \\\frac{\mu\nu(a-b)}{1-\mu\nu}\endcases,\quad G_1=\cases \frac 1{1-\mu\nu}\\ \frac\mu{1-\mu\nu}\endcases,\quad G_2=\cases \frac\nu{1-\mu\nu}\\ \frac 1{1-\mu\nu}\endcases,\\
\lim\limits_{\ze\to a}\left(m'(\zeta)\bmatrix f_1(a)\\f_2(a)\endbmatrix\right)=\bmatrix \frac{-\mu\nu}{1-\mu\nu}\\ \frac{\mu^2\nu}{1-\mu\nu}\endbmatrix,\quad
\lim\limits_{\ze\to b}\left(m'(\zeta)\bmatrix f_1(b)\\ g_1(b)\endbmatrix\right)=\bmatrix\frac{\mu\nu^2}{1-\mu\nu}\\ \frac{-\mu\nu}{1-\mu\nu}\endbmatrix.
\endgathered
$$
By Theorem 3.5, this implies that
$$
K=\frac L{1+L}=\frac 1{1-\mu\nu}\bmatrix  -\mu\nu &\mu\\ \nu &-\mu\nu \endbmatrix
$$
which is immediately verified directly. Note that the condition $1-\mu\nu\ne 0$ is equivalent to the invertibility of $(1+L)$.
\endexample

In what follows we will demonstrate how to use Theorems 3.3 and 3.5 to compute correlation kernels of determinantal point processes arising in concrete representation theoretic models. 

\head 4. Harmonic analysis on $S(\infty)$: Whittaker kernel \endhead

As is explained in \cite{BO2}, see also \cite{Ol2}, the problem of decomposing generalized regular representations of the infinite symmetric group $S(\infty)$ on irreducible ones reduces to computing correlation kernels of certain determinantal point processes.

Specifically, consider a determinantal point process on $\Z'=\Z+\frac 12$ constructed using Example 1.2 with the L-kernel given by
$$
L(x,y)=\cases 0,&xy>0,\\
\dfrac{|z(z+1)_{x-\frac 12}(-z+1)_{-y-\frac 12}|\,\xi^{\frac{x-y}2}}{(x-\frac 12)!(-y-\frac 12)!(x-y)}\,,&x>0,\,y<0,\\
\dfrac{|z(-z+1)_{-x-\frac 12}(z+1)_{y-\frac 12}|\,\xi^{\frac{-x+y}2}}{(-x-\frac 12)!(y-\frac 12)!(x-y)}\,,&x<0,\,y>0.
\endcases
$$
Here $z\in\C\setminus\Z$ and $\xi\in(0,1)$ are parameters. The symbol $(a)_k$ stands for $a(a+1)\cdots(a+k-1)=\Gamma(a+k)/\Gamma(a)$.

Note that as $|z|\to\infty$, $\xi\to 0$, and $|z|^2\xi\to\theta$, this kernel converges to the L-kernel of Example 2.2.

The problem consists in computing $K=L(1+L)^{-1}$ and taking the scaling limit
$$
\Cal K (x,y)=\lim_{\xi\to 1}(1-\xi)^{-1}\cdot K\left([(1-\xi)^{-1}x]+\tfrac12,\,[(1-\xi)^{-1}y]+\tfrac 12\right).
$$

This problem has been solved in \cite{BO2}. However, we did not provide a {\it derivation} of the formula for the kernel $K$ there, we just {\it verified} the equality $K=L(1+L)^{-1}$.

The goal of this section is to provide a derivation of the kernel 
$\Cal K(x,y)$ bypassing the computation of $K(x,y)$.

Observe that there exists a limit
$$
\Cal L(x,y)=\lim_{\xi\to 1} (1-\xi)^{-1}\cdot
L\left([(1-\xi)^{-1}x]+\tfrac12,\,[(1-\xi)^{-1}y]+\tfrac 12\right).
$$
Indeed, for $a=[(1-\xi)^{-1}x]$, $b=[(1-\xi)^{-1}y]$,
$$
\gathered
\frac{(z+1)_a}{\Gamma(a+1)}=\frac{\Gamma(z+1+a)}{\Gamma(z+1)\Gamma(a+1)}\sim \frac{a^z}{\Gamma(z+1)}, \quad x>0,\\
\frac{(-z+1)_b}{\Gamma(b+1)}=\frac{\Gamma(-z+1+|b|)}{\Gamma(-z+1)\Gamma(|b|+1)}\sim \frac{|b|^{-z}}{\Gamma(-z+1)}, \quad y<0,\\
\xi^{\frac a2}\sim (1-(1-\xi))^{\frac x{2(1-\xi)}}\sim e^{\frac x2},\quad
\xi^{-\frac b2}\sim e^{-\frac y2},
\endgathered
$$
and we get
$$
\Cal L(x,y)=\cases 0,&xy>0,\\
\dfrac{|\sin\pi z|}\pi \, \dfrac{(x/|y|)^{\Re z}e^{\frac{-x+y}2}}{x-y}\,,&
x>0,\,y<0,\\
\dfrac{|\sin\pi z|}\pi \, \dfrac{(y/|x|)^{\Re z}e^{\frac{x-y}2}}{x-y}\,,&
x<0,\,y>0.
\endcases
$$

It is natural to assume that $\Cal K=\Cal L(1+\Cal L)^{-1}$. It turns out that this relation holds whenever $\Cal L$ defines a bounded operator in $L^2(\R\setminus \{0\})$, which happens when $|\Re z|<\frac 12$, see \cite{Ol1}. Our goal is to derive $\Cal K$ using the relation $\Cal L=\Cal L(1+\Cal L)^{-1}$. 

It is easily seen that $\Cal L$ is an integrable operator; we can take
$$
\gathered
f_1(x)=g_2(x)=\cases \dfrac{|z|^\frac 12}{|\Gamma(z+1)|}\, x^{\Re z}e^{-\frac x2},&x>0,\\0,&x<0,\endcases\\
f_2(x)=g_1(x)=\cases 0,&x>0,\\ \dfrac{|z|^\frac 12}{|\Gamma(-z+1)|}\, |x|^{-\Re z}e^{\frac x2},&x<0.\endcases
\endgathered
$$
Note that $\Cal L(y,x)=-\Cal L(x,y)$ which means that $(1+\Cal L)$ is invertible (provided that $\Cal L$ is bounded).

The RHP of Theorem 3.3 then has the jump matrix
$$
v(x)=\cases \bmatrix 1& 2i\,|\sin\pi z\,|\,x^{2\Re z}e^{-x}\\0&1\endbmatrix,&x>0\\
\bmatrix 1&0\\ 2i\,|\sin\pi z\,|\, |x|^{-2\Re z}e^x&1\endbmatrix, &x<0.
\endcases
$$

The key property of this RHP which will allow us to solve it, is that it can be reduced to a problem with a piecewise constant jump matrix.

Let $m$ be the solution of the normalized RHP $(\R\setminus\{0\},v)$. Set
$$
\Psi(\zeta)=m(\ze)\bmatrix \ze^{\Re z}e^{-\frac \ze 2}&0\\0&\ze^{-\Re z}e^{\frac \ze 2}\endbmatrix,\qquad \ze\notin \R.
$$
Then the jump relation $m_+=m_-v$ takes the form
$$
\Psi_+(x)=\Psi_-(x)
\bmatrix x^{-\Re z}e^{\frac x 2}&0\\0&x^{\Re z}e^{-\frac x 2}\endbmatrix_-v(x)
\bmatrix x^{\Re z}e^{-\frac x 2}&0\\0&x^{-\Re z}e^{\frac x 2}\endbmatrix_+,
$$
and a direct computation shows that the jump matrix for $\Psi$ takes the form
$$
\cases \bmatrix 1& 2i\,|\sin\pi z\,|\\0&1\endbmatrix,&x>0\\
\bmatrix e^{2\pi i \,\Re z}&0\\ 2i\,|\sin\pi z\,|&e^{-2\pi i\, \Re z} \endbmatrix, &x<0.
\endcases
$$
Let us first find a solution of this RHP without imposing any asymptotic conditions at infinity. We will denote it by $\Psi^0$. Set
$$
\wt \Psi^0{\ze}=\cases \Psi^0(\ze),&\Im\ze>0,\\
\Psi^0(\ze)\bmatrix 1&2i\,|\sin\pi z|\\0&1\endbmatrix,&\Im\ze<0.
\endcases
$$
Then $\wt\Psi^0$ has no jump across $\R_+$, and the jump matrix $(\wt\Psi^0_-)^{-1}\wt\Psi^0_+$ on $\R_-$ has the form
$$
\multline
\bmatrix 1&-2i\,|\sin\pi z|\\0&1\endbmatrix
\bmatrix e^{2\pi i \,\Re z}&0\\ 2i\,|\sin\pi z\,|&e^{-2\pi i\, \Re z} \endbmatrix\\=
\bmatrix e^{2\pi i \,\Re z}+4|\sin\pi z\,|^2& -2i\,|\sin\pi z|e^{-2\pi i \,\Re z}
\\ 2i\,|\sin\pi z\,|&e^{-2\pi i\, \Re z} \endbmatrix.
\endmultline
$$
The determinant of this matrix is equal to 1, and the trace is equal to
$2\cos(2\pi\Re z)+4|\sin\pi z\,|^2=e^{\pi\Im z}+e^{-\pi\Im z}.$ Thus, if $z\notin \R$, there exists a nondegenerate $U$ such that 
$$
\bmatrix e^{2\pi i \,\Re z}+4|\sin\pi z\,|^2& -2i\,|\sin\pi z|e^{-2\pi i \,\Re z}
\\ 2i\,|\sin\pi z\,|&e^{-2\pi i\, \Re z} \endbmatrix=
U^{-1}\bmatrix e^{\pi\Im z}&0\\0& e^{-\pi\Im z}\endbmatrix U.
$$ 
This means that the matrix $\wt\Psi^0 U^{-1}$ has jump $\bmatrix e^{\pi\Im z}&0\\0& e^{-\pi\Im z}\endbmatrix$ across $\R_-$. Note that the matrix
$\bmatrix \ze^{-i\Im z}&0\\0& \ze^{i\Im z}\endbmatrix$ satisfies the same jump relation. Hence,
$$
\Psi^0(\ze)=\cases \bmatrix \ze^{-i\Im z}&0\\0& \ze^{i\Im z}\endbmatrix U,&\Im\ze>0,\\
\bmatrix \ze^{-i\Im z}&0\\0& \ze^{i\Im z}\endbmatrix\bmatrix 1&-2i\,|\sin\pi z\,|\\0&1\endbmatrix U,&
\Im \ze<0,
\endcases
$$
is a solution of our RHP for $\Psi$. It follows that $\Psi(\Psi^0)^{-1}$ has no jump across $\R$ and this implies, modulo some technicalities, that $\Psi=H\Psi^0$ where $H$ is entire.

Now we describe the crucial step. Since the jump matrix for $\Psi$ is piecewise constant, $\Psi'=\frac{d\Psi}{d\ze}$ satisfies the same jump condition as $\Psi$, and hence $\Psi'\Psi^{-1}$ is meromorphic in $\C$ with a possible pole at $\ze=0$. On the other hand we have
$$
\gathered
\Psi'\Psi^{-1}=H'H^{-1}+H(\Psi^0)'(\Psi^0)^{-1}H^{-1}\\ =H'H^{-1}
-\frac 1\zeta\, i\,\Im z\, H\bmatrix 1&0\\0&-1\endbmatrix H^{-1}=
\text{entire function}+\frac A\zeta,
\endgathered
\tag 4.1
$$
where $A$ has eigenvalues $\pm i\,\Im z$.

Let us recall now that $m$ solves the {\it normalized} RHP, which means that $m(\ze)\sim I$ as $|z|\to\infty$. An additional argument shows that
$$
m(\ze)=I+{m^{(1)}}\ze^{-1}+O\left({|\ze|^{-2}}\right),\qquad |\ze|\to\infty,
$$
with a constant matrix $m^{(1)}$. Thus,
$$
\Psi(\ze)=\left(I+{m^{(1)}}\ze^{-1}+O\left({|\ze|^{-2}}\right)\right)
\bmatrix \ze^{\Re z}e^{-\frac \ze 2}&0\\0&\ze^{-\Re z}e^{\frac \ze 2}\endbmatrix
$$
and
$$
\Psi'(\ze)\Psi^{-1}(\ze)=-\frac 12\bmatrix 1&0\\0&-1\endbmatrix
+\frac 1\ze\bmatrix \Re z&-m^{(1)}_{12}\\m_{21}^{(1)}&-\Re z\endbmatrix
+O(|\ze|^{-2}).
$$
Comparing this relation with \tht{4.1} we conclude that
$$
\Psi'(\ze)=\left(-\frac 12\bmatrix 1&0\\0&-1\endbmatrix
+\frac 1\ze\bmatrix \Re z&-m^{(1)}_{12}\\m_{21}^{(1)}&-\Re z\endbmatrix\right)\Psi(\ze)
$$
with $m^{(1)}_{12}m^{(1)}_{21}=(\Re z)^2+(\Im \ze)^2=|z|^2$. This 1st order linear matrix differential equation leads to 2nd order linear differential equations on the matrix elements on $\Psi$, for example
$$
\zeta\Psi_{11}''+\Psi_{11}'=\left(-\frac 12-\frac{|z|^2}{\ze}+\frac 1\ze\left(\Re z-\frac \ze 2\right)^2\right)\Psi_{11}.
$$
Using these differential equations and the asymptotics of $\Psi$ at infinity, it is easy to express $\Psi$ in terms of the confluent hypergeometric function or the Whittaker function, see \cite{Er, 6.9} for definitions.
In terms of the Whittaker function $W_{\kappa,\mu}$ the final formula for $\Psi$ has the form
$$
\Psi(\ze)= \bmatrix \ze^{-\frac{1}2}W_{\Re z+\frac{1}2,i\Im z}(\ze)&
|z|\,(-\ze)^{-\frac{1}2}W_{-\Re z-\frac{1}2,i\Im z}(-\ze)\\
-|z|\,\ze^{-\frac{1}2}W_{\Re z-\frac{1}2,i\Im z}(\ze)&(-\ze)^{-\frac{1}2}W_{-\Re z+\frac{1}2,i\Im z}(-\ze)\endbmatrix.
$$
It is not hard to show that $\det \Psi\equiv 1$, and
$$
\Psi^{-t}(\ze)= \bmatrix (-\ze)^{-\frac{1}2}W_{-\Re z+\frac{1}2,i\Im z}(-\ze)&|z|\,\ze^{-\frac{1}2}W_{\Re z-\frac{1}2,i\Im z}(\ze)
\\
-|z|\,(-\ze)^{-\frac{1}2}W_{-\Re z-\frac{1}2,i\Im z}(-\ze)&\ze^{-\frac{1}2}W_{\Re z+\frac{1}2,i\Im z}(\ze)\endbmatrix.
$$
Then Theorem 3.3 implies
$$
\gathered
F_1(x)=\cases \frac{|z|^\frac 12}{|\Gamma(z+1)|}\,\Psi_{11}(x),&x>0,\\
\frac{|z|^\frac 12}{|\Gamma(-z+1)|}\,\Psi_{12}(x),&x<0,
\endcases
\quad
F_2(x)=\cases \frac{|z|^\frac 12}{|\Gamma(z+1)|}\,\Psi_{21}(x),&x>0,\\
\frac{|z|^\frac 12}{|\Gamma(-z+1)|}\,\Psi_{22}(x),&x<0,
\endcases\\
G_1(x)=\cases -\frac{|z|^\frac 12}{|\Gamma(z+1)|}\,\Psi_{21}(x),&x>0,\\
\frac{|z|^\frac 12}{|\Gamma(-z+1)|}\,\Psi_{22}(x),&x<0,
\endcases
\quad
G_2(x)=\cases \frac{|z|^\frac 12}{|\Gamma(z+1)|}\,\Psi_{11}(x),&x>0,\\
-\frac{|z|^\frac 12}{|\Gamma(-z+1)|}\,\Psi_{12}(x),&x<0,
\endcases
\endgathered
$$
and
$$
\Cal K(x,y)=\frac{F_1(x)G_1(y)+F_2(x)G_2(y)}{x-y}\,,\quad x,y\in\R\setminus\{0\}.
$$
This kernel is called the {\it Whittaker kernel}, see \cite{BO1}, \cite{B1}. 

\head 5. Poissonized Plancherel measure: discrete Bessel kernel
\endhead

We now return to the situation described in Example 2.2. Our goal is to compute the correlation kernel $K=L(1+L)^{-1}$. The exposition below follows \cite{B2, \S7}.

According to Theorem 3.5, we have to find the unique solution of the normalized DRHP $(\Z',w)$ with
$$
w(x)=\cases \bmatrix 0&-\frac{\theta^x}{\left(x-\frac 12\right)!^2}\\0&0\endbmatrix,
& x\in\Z'_+,\\
\bmatrix 0&0\\ -\frac{\theta^{-x}}{\left(-x-\frac 12\right)!^2}&0\endbmatrix,&
x\in\Z'_-.
\endcases
$$
Note that the kernel $L$ is skew-symmetric, which means that $(1+L)$ is invertible. If we denote by $m$ the solution of this DRHP then
$$
m(\ze)=I+\bmatrix \alpha&\beta\\ \gamma&\delta\endbmatrix\zeta^{-1}+
O(|\zeta^{-2}|),\qquad |\ze|\to\infty,
$$
with constant $\alpha,\,\beta,\,\gamma,\,\delta$.
The symmetry of the problem with respect to 
$$
\ze\leftrightarrow -\ze, \quad \bmatrix m_{11}(\ze) & m_{12}(\ze)\\ m_{21}(\ze) &m_{22}(\ze)\endbmatrix\longleftrightarrow
\bmatrix m_{22}(-\ze) &- m_{21}(-\ze)\\- m_{12}(-\ze)& m_{11}(-\ze)\endbmatrix,
$$
implies that $\gamma=\beta$ and $\delta=-\alpha$. 

Denote $\eta=\sqrt{\theta}$ and set  
$$
n(\ze)=m(\ze)\bmatrix \eta^\ze&0\\0&\eta^{-\ze}\endbmatrix.
$$
Then $n(\ze)$ solves a DRHP with the jump matrix
$$
\cases \bmatrix 0&-\frac{1}{\left(x-\frac 12\right)!^2}\\0&0\endbmatrix,
& x\in\Z'_+,\\
\bmatrix 0&0\\ -\frac{1}{\left(-x-\frac 12\right)!^2}&0\endbmatrix,&
x\in\Z'_-.
\endcases
$$
Note that this jump matrix does not depend on $\eta$. This means that
$\frac{\partial n}{\partial\eta}$ has the same jump matrix, and hence
the matrix $\frac{\partial n}{\partial\eta}\,n^{-1}$ is entire. Computing the asymptotics as $\ze\to\infty$, we obtain
$$
\frac{\partial n(\ze)}{\partial\eta}\,n^{-1}(\ze)=\bmatrix
\zeta&-2\beta\\2\beta&-\zeta\endbmatrix+O(|z|^{-1}).
$$
By Liouville's theorem, the remainder term must vanish, and thus
$$
\frac{\partial n(\ze)}{\partial\eta}=\bmatrix
\zeta&-2\beta\\2\beta&-\zeta\endbmatrix n(\ze).
\tag 5.1
$$
This yields 2nd order linear differential equations on the matrix elements of $n$ which involve, however, an unknown function $\beta=\beta(\eta)$.

In order to determine $\beta$ we need to make one more step. Set
$$
p(\ze)=n(\ze)\bmatrix \frac 1{\Gamma(\ze+\frac 12)}&0\\ 0 &\frac 1{\Gamma(-\ze+\frac 12)}\endbmatrix.
$$
It is immediately verified that the fact that $n$ solves the corresponding DRHP is equivalent to $p$ being entire and satisfying the condition
$$
p(x)=(-1)^{x-\frac 12}p(x)\bmatrix 0&1\\1&0\endbmatrix,\qquad x\in\Z'.
\tag 5.2
$$
The key property of this relation is that it depends on $x$ in an insubstantial way. This allows us to do the following trick which should be viewed as a substitute of the differentiation with respect to $x$. 
Set
$$
\wt p(\ze)=\bmatrix p_{11}(\ze+1)&-p_{12}(\ze+1)\\-p_{21}(\ze-1)&p_{22}(\ze-1)\endbmatrix.
$$
Then $\wt p$ satisfies the same condition \tht{5.2} as $p$ does, and hence 
$$
\wt n(\ze)=\wt p(\ze) \bmatrix {\Gamma(\ze+\frac 12)}&0\\ 0 &{\Gamma(-\ze+\frac 12)}\endbmatrix
$$
satisfies the same DRHP as $n$ does. Thus, $\wt n\,n^{-1}$ is entire. Working out the asymptotics as $\ze\to \infty$, we obtain
$$
\wt n(\zeta)n^{-1}(\ze)=\bmatrix 0&\frac \beta\eta\\-\frac\beta\eta&0\endbmatrix
+O(|z|^{-1}).
$$
Liouville's theorem implies that the remainder vanishes. Hence, 
$$
\wt p_{11}(\ze)=p_{11}(\ze+1)=\frac \beta\eta \,p_{21}(\ze),\quad
\wt p_{21}(\ze)=-p_{21}(\ze-1)=-\frac\beta\eta \, p_{11}(\ze).
$$
This implies that $(\beta/\eta)^2=1$. Both cases $\beta=\pm\eta$ lead, via \tht{5.1}, to linear 2nd order differential equations on the matrix elements of $n$ or matrix elements of $p$. For example, $\beta=-\eta$ yields
$$
\left(\frac{\partial^2}{\partial\eta^2}-\frac{\ze(\ze-1)}{\eta^2}+4\right)p_{11}(\ze)=0,\quad \left(\frac{\partial^2}{\partial\eta^2}-\frac{\ze(\ze+1)}{\eta^2}+4\right)p_{21}(\ze)=0.
$$
General solutions of these equations can be written in terms of Bessel functions, and matching the asymptotics at infinity we obtain for $\beta=-\eta$
$$
p(\ze)=\sqrt{\eta}\bmatrix J_{\ze-\frac 12}(2\eta)&J_{-\ze+\frac 12}(2\eta)\\
-J_{\ze+\frac 12}(2\eta)&J_{-\ze-\frac 12}(2\eta)\endbmatrix,
\tag 5.3
$$
and for $\beta=\eta$
$$
\widehat p(\ze)=\sqrt{\eta}\bmatrix J_{\ze-\frac 12}(2\eta)&-J_{-\ze+\frac 12}(2\eta)\\
J_{\ze+\frac 12}(2\eta)&J_{-\ze-\frac 12}(2\eta)\endbmatrix.
$$
Here $J_{\nu}(u)$ is the Bessel function, see \cite{Er, 7.2} for the definition.

Using the well--known relation $J_{-n}=(-1)^n J_n$ we immediately see that $p(\ze)$ given by \tht{5.3} satisfies \tht{5.2}, while $\widehat p(\ze)$ does not. In fact $\widehat p(\ze)$ satisfies \tht{5.2} with 
$(-1)^{x-\frac 12}$ replaced with $(-1)^{x+\frac 12}$. This means that
$$
m(\ze)=p(\ze)\bmatrix {\eta^{-\ze}\Gamma(\ze+\frac 12)}&0\\ 0 &{\eta^{\ze}\Gamma(-\ze+\frac 12)}\endbmatrix
$$
solves the initial DRHP, and by Theorem 3.5 we obtain
$$
\gathered
F_1(x)=\cases p_{11}(x),&x>0,\\
p_{12}(x),&x<0,
\endcases
\quad
F_2(x)=\cases p_{21}(x),&x>0,\\
p_{22}(x),&x<0,
\endcases\\
G_1(x)=\cases -p_{21}(x),&x>0,\\
p_{22}(x),&x<0,
\endcases
\quad
G_2(x)=\cases p_{11}(x),&x>0,\\
-p_{12}(x),&x<0,
\endcases
\endgathered
$$
and
$$
K(x,y)=\frac{F_1(x)G_1(y)+F_2(x)G_2(y)}{x-y}\,,\quad x,y\in\Z'.
$$
The diagonal values $K(x,x)$ are determined by the L'Hospital rule:
$$
K(x,x)=F_1'(x)G_1(x)+F_2'(x)G_2(x).
$$

This is the discrete Bessel kernel obtained in \cite{BOO}. The restriction of $K(x,y)$ to $\Z'_+\times \Z'_+$ was independently derived in \cite{J}.

It is worth noting that the matrix $\widehat p$ also has an important meaning. In fact, if we define a kernel $\widehat K$ using the formulas above with $p$ replaced by $\widehat p$ then $\widehat K=L(L-1)^{-1}$.

\vskip 1 true cm

School of Mathematics, Institute for Advanced
Study, Einstein Drive, Princeton NJ 08540, U.S.A.

E-mail address:
{\tt borodine\@math.upenn.edu}

\enddocument